\documentclass{article}

\usepackage{amsmath}
\usepackage{amssymb}

\newtheorem{thm}{Theorem}
\newtheorem{conj}{Conjecture}
\newtheorem{lem}{Lemma}

\newtheorem{cor}{Corollary}

\title{Squarefull Numbers in Arithmetic Progression II}
\author{Tsz Ho Chan}
\date{}

\begin{document}
\maketitle

\begin{abstract}
In this paper, we improve the error term in a previous paper on an asymptotic formula for the number of squarefull numbers in an arithmetic progression.
\end{abstract}

\section{Introduction and Main Results}
A positive integer $n$ is squarefull if $n = 1$ or for every prime
$p$ dividing $n$, $p^2$ also divides $n$. They are also known as
powerful numbers. In [\ref{CT}], Professor Tsang and the author studied squarefull numbers in arithmetic progression and obtained
\begin{thm} \label{oldthm1}
For $(l,q) = 1$,
\[
T_q(x; l) := \mathop{\mathop{\sum_{n \le x}}_{n \text{ squarefull}}}_{n
\equiv l \; (\bmod q)} 1 = \frac{C_q(l)}{\zeta(3)} \prod_{p | q}
\Bigl( \frac{p^3}{p^3 - 1} \Bigr) \frac{x^{1/2}}{q} + \frac{D_q(l)}{\zeta(2)}
\prod_{p | q} \Bigl( \frac{p^2}{p^2 - 1} \Bigr) \frac{x^{1/3}}{q} + O_\epsilon (\min(x^{1/6} q^{1/3 + \epsilon},
x^{1/5} q^{\epsilon})).
\]
\end{thm}
Here, for any modulus $q$, we define
\[
N_2(n ; q) := \# \{ x \; (\bmod \; q) : x^2 \equiv n (\bmod \; q),
(x,q) = 1 \}
\]
and
\[
N_3(n ; q) := \# \{ x \; (\bmod \; q) : x^3 \equiv n (\bmod \; q),
(x,q) = 1 \}
\]
where $\# S$ stands for the number of elements in a set $S$. Also 
\[
C_q(l) := \sum_{b = 1}^{\infty} \frac{N_2(l b; q)}{b^{3/2}}
\]
and
\[
D_q(l) := -2 \frac{\phi(q)}{q} + \frac{2}{3}
\int_{1}^{\infty} \frac{\sum_{a \le u} N_3(l a; q)
- \frac{\phi(q)}{q} u}{u^{5/3}} du.
\]
Note that $N_2(n a^2; q) = N_2(n; q)$ and $N_3(n a^3; q) = N_3(n; q)$ for
any $(a,q) = 1$, and $N_2(n; q) = 0 = N_3(n; q)$ if $(n,q) > 1$. Also note that
\begin{equation} \label{full}
\sum_{n \le q} N_2(n; q) = \phi(q) \text{ and } \sum_{n \le q} N_3(n,q) = \phi(q).
\end{equation}

More recently, Srichan [\ref{S}] used the method of exponent pairs and obtained the same result with a slightly bigger error term.
However, he expressed the main terms using Dirichlet characters and Dirichlet $L$-functions. He also studied the case where $l$ and $q$ are not
relatively prime as well as cubefull numbers in arithmetic progressions.

In this paper, we improve the error term further by Kloosterman-type exponential sum method and obtain
\begin{thm} \label{finalthm}
For $(l,q) = 1$,
\[
T_q(x; l) = \frac{C_q(l)}{\zeta(3)} \prod_{p | q}
\Bigl( \frac{p^3}{p^3 - 1} \Bigr) \frac{x^{1/2}}{q} + \frac{D_q(l)}{\zeta(2)}
\prod_{p | q} \Bigl( \frac{p^2}{p^2 - 1} \Bigr) \frac{x^{1/3}}{q} + O_\epsilon \Bigl( \Bigl( x^{1/6} q^{1/12} + \frac{x^{1/5}}{q^{1/5}} \Bigr) q^\epsilon \Bigr).
\]
\end{thm}
As an intermediate step, we also prove
\begin{thm} \label{interthm}
For $(l,q) = 1$ and $0 < \mu < 1/5$,
\[
S_q(X;l) := \mathop{\sum_{a^2 b^3 \le X}}_{a^2 b^3 \equiv l \; (\bmod q)} 1 = C_q(l) \frac{X^{1/2}}{q} + D_q(l) \frac{X^{1/3}}{q} + O_\epsilon \Bigl(\frac{X^{1/2 - 3\mu/2}}{q^{1/2 - \epsilon}} + \min(X^\mu, q^{1/2}) q^\epsilon \Bigr).
\]
\end{thm}
As an immediate corollary, we have
\begin{cor}
For $(l,q) = 1$,
\[
T_q(x; l) \ll_\epsilon \frac{x^{1/2 + \epsilon}}{q} \; \;
\text{ for } \; \; q \le x^{1/3}.
\]
\end{cor}
However, one suspects the following is true.
\begin{conj}
For $(l,q) = 1$,
\[
T_q(x; l) \ll_\epsilon \frac{x^{1/2 + \epsilon}}{q} \; \;
\text{ for } \; \; q \le x^{1/2}.
\]
\end{conj}
To prove Theorem \ref{interthm}, we need a result on the number of points of the curve $x^2 y^3 \equiv l (\bmod q)$ in a box.
\begin{thm} \label{firstthm}
For $(l,q) = 1$,
\[
N_{A,B}(K,L) =: \# \{(a,b) : a^2 b^3 \equiv l (\bmod q), A < a \le A+K, B < b \le B+L \} = \frac{\phi(q) K L}{q^2} + O\Bigl(\frac{K+L}{q^{1/2-\epsilon}} + q^{1/2+\epsilon}\Bigr).
\]
\end{thm}

The paper is organized as follows. First, we state and prove some lemmas which are needed later. Then we prove Theorem \ref{firstthm} and use it to prove Theorem \ref{interthm}. Finally, we prove Theorem \ref{finalthm}.

\bigskip

{\bf Some Notations} Throughout the paper, $p$ stands for a prime.
$p | q$ means that $p$ divides $q$, and $p^a || q$ means that $p^a |
q$ but $p^{a+1} \nmid q$. The symbol $\overline{a}$ stands for the
multiplicative inverse of $a \pmod q$ (i.e. $a \overline{a} \equiv 1
\pmod q$). The notations $f(x) = O(g(x))$, $f(x) \ll g(x)$ and $g(x)
\gg f(x)$ are all equivalent to $|f(x)| \leq C g(x)$ for some
constant $C > 0$. Finally $f(x) = O_\lambda(g(x))$, $f(x)
\ll_\lambda g(x)$ or $g(x) \gg_\lambda f(x)$ mean that the implicit
constant $C$ may depend on $\lambda$.
\section{Lemmas}
\begin{lem} \label{lemgauss}
Let $\chi$ be a character $(\bmod q)$ and $\chi^* (\bmod q^*)$ be the primitive character that induces $\chi$. Then
\[
\Big| \sum_{x = 1}^{q} \chi(x) e\Bigl(\frac{a x}{q}\Bigr) \Big| \le \sigma(a, q/q^*) \sqrt{q^*}.
\]
\end{lem}

Proof: See equation (12.48) on page 324 of [\ref{IK}].

\begin{lem} \label{lem1}
For $(l,q)=1$, $\alpha \le 1$,
\[
S =: \sum_{n = 1}^{q-1} (n,q) \Big| \sum_{b = B+1}^{B+L} e\Bigl(-\frac{n b}{q}\Bigr) \Big| \ll d(q) q \log q.
\]
\end{lem}

Proof: Let $\|x\|$ be the distance from $x$ to the nearest integer. Then we have
\begin{align*}
S =& \sum_{d | q} d \sum_{n' = 1}^{q/d - 1} \Big| \sum_{b = B+1}^{B+L} e\Bigl(-\frac{n' b}{q/d}\Bigr) \Big| \\
&\ll \sum_{d | q} d \sum_{n' = 1}^{q/d - 1} \min\Bigl(L, \frac{1}{\|\frac{n'}{q/d}\|}\Bigr) \\
&\ll \sum_{d | q} d \Bigl[ \sum_{n' \le q/(Ld)} L + \sum_{q/(Ld) < n' \le q/(2d)} \frac{q/d}{n'} \Bigr] \ll \sum_{d | q} d \frac{q}{d} \log q \ll d(q) q \log q.
\end{align*}

\begin{lem} \label{lem2}
For $(l,q)=1$,
\[
\sum_{b \le u} N_2(lb; q) - \frac{\phi(q)}{q} u = O_\epsilon(q^{1/2 + \epsilon}).
\]
\[
\sum_{a \le u} N_3(la; q) - \frac{\phi(q)}{q} u = O_\epsilon(q^{1/2 + \epsilon}).
\]
\end{lem}

Proof: By (\ref{full}), we can assume that $u \le q$. By orthogonal property of additive characters, we have
\begin{align*}
\sum_{b \le u} N_2(lb; q) =& \sum_{b \le u} \frac{1}{q} \sum_{n = 1}^{q} \mathop{\sum' \nolimits}_{x = 1}^{q} e\Bigl(\frac{n(x^2 - lb)}{q}\Bigr) \\
=& \frac{\phi(q)}{q} (u + O(1)) + \frac{1}{q} \sum_{n = 1}^{q-1} \sum_{b \le u} e\Bigl(\frac{-nlb}{q}\Bigr) \mathop{\sum' \nolimits}_{x = 1}^{q} e\Bigl(\frac{n x^2}{q}\Bigr) \\
=& \frac{\phi(q)}{q} (u + O(1)) + \frac{1}{q} \sum_{n = 1}^{q-1} \sum_{b \le u} e\Bigl(\frac{-nlb}{q}\Bigr) \sum_{y (\bmod q)} e\Bigl(\frac{n y}{q} \Bigr) \sum_{\chi \in G_2} \chi(y)
\end{align*}
where $G_2$ is the set of all characters $\chi (\bmod q)$ such that $\chi^2 = \chi_0$, the principal character. Interchanging the summations, we have
\begin{align*}
\sum_{b \le u} N_2(lb; q) =& \frac{\phi(q)}{q} (u + O(1)) + \frac{1}{q} \sum_{\chi \in G_2} \sum_{n = 1}^{q-1} \sum_{b \le u} e\Bigl(\frac{-nlb}{q}\Bigr) \sum_{y (\bmod q)} \chi(y) e\Bigl(\frac{n y}{q} \Bigr) \\
=& \frac{\phi(q)}{q} (u + O(1)) + O\Bigl(\frac{1}{q} \sum_{\chi \in G_2} \sum_{n = 1}^{q-1} \Big| \sum_{b \le u} e\Bigl(\frac{-nlb}{q}\Bigr)\Big| (n,q) d((n,q)) \sqrt{q} \Bigr) \\
=& \frac{\phi(q)}{q} (u + O(1)) + O\Bigl(\frac{d(q)}{\sqrt{q}} |G_2| \sum_{n = 1}^{q-1} \Big| \sum_{b \le u} e\Bigl(\frac{-nlb}{q}\Bigr)\Big| (n,q) \Bigr) \\
=& \frac{\phi(q)}{q} u + O_\epsilon(q^{1/2 + \epsilon})
\end{align*}
by Lemma \ref{lemgauss}, Lemma \ref{lem1}, and $d(q) \ll_\epsilon q^\epsilon$, $N_2(n; q) \ll_\epsilon q^{\epsilon}$ (See Nagell [\ref{N}]) and the fact that the group of characters $(\bmod q)$ is isomorphic to the group of reduced residues $\mathbb{Z}/q\mathbb{Z}^{*}$. The proof of the result for $N_3(l a; q)$ is almost identical except that one should use $G_3$, the set of all characters $\chi (\bmod q)$ such that $\chi^3 = \chi_0$, instead of $G_2$.

\bigskip

Below we are going to study the exponential sum
\[
S(a,b;q) := \mathop{\sum' \nolimits}_{n = 1}^{q} e\Bigl(\frac{a n^2 + b \overline{n}^3}{q}\Bigr).
\]
Suppose $q = rs$ with $(r,s)=1$. By the ``reciprocity" formula
\[
\frac{\overline{s}}{r} + \frac{\overline{r}}{s} = \frac{1}{q} (\bmod 1)
\]
where $s \overline{s} \equiv 1 (\bmod r)$ and $r \overline{r} \equiv 1 (\bmod s)$,
\[
e\Bigl(\frac{a n^2 + b \overline{n}^3}{q}\Bigr) = e\Bigl(\frac{a \overline{s} n^2 + b \overline{s} \, \overline{n}^3}{r}\Bigr) e\Bigl(\frac{a \overline{r} n^2 + b \overline{r} \, \overline{n}^3}{s}\Bigr)
\]
and hence with $n = s x + r y$ with $1 \le x \le r$, $1 \le y \le s$, $(x,r)=1$ and $(y,s)=1$
\begin{equation} \label{chain}
S(a,b;q) = \mathop{\sum' \nolimits}_{x = 1}^{r} \mathop{\sum' \nolimits}_{y = 1}^{s} e\Bigl(\frac{a \overline{s} x^2 + b \overline{s} \, \overline{x}^3}{r}\Bigr) e\Bigl(\frac{a \overline{r} y^2 + b \overline{r} \, \overline{y}^3}{s}\Bigr) = S(a \overline{s}, b \overline{s}; r) S(a \overline{r}, b \overline{r}; s).
\end{equation}
Therefore the study of $S(a,b;q)$ reduces to that of prime power moduli. When $q = p$ is prime, we have by Theorem 5 of [\ref{Bom}] with the rational function $R(x) = \frac{a x^5 + b}{x^3}$ that
\begin{equation} \label{exp-prime}
|S(a,b;p)| \le 2 (a,b,p)^{1/2} p^{1/2}.
\end{equation}
Note: When $(a,b,p)=1$, one needs to check that $R(x) \neq h(x)^p - h(x)$ for $h(x) \in \overline{\mathbb{Z}/p\mathbb{Z}}[x]$. Suppose not, then
\begin{equation} \label{exp-check}
\frac{a x^5 + b}{x^3} = \frac{f(x)^p}{g(x)^p} - \frac{f(x)}{g(x)} \text{ or } g(x)^p (a x^5 + b) = x^3 (f(x)^p - g(x)^{p-1} f(x))
\end{equation}
with $(f(x),g(x))=1$. This forces $g(x)^p | x^3$. If $p > 3$, then this forces $g(x)$ to be a constant. Then (\ref{exp-check}) is impossible by degree consideration. If $p \le 3$, one can obtain (\ref{exp-prime}) by bounding the exponential sum trivially.

\bigskip

For prime power moduli $q = p^\beta$ with $\beta \ge 2$, we may assume that $(a,b,p)=1$ for otherwise we can reduce the modulus to a lower exponent. If $\beta = 2 \alpha$ with $\alpha \ge 1$, by Lemma 12.2 in [\ref{IK}], we have
\begin{equation} \label{exp-evenppower}
S(a,b;p^{2 \alpha}) = p^\alpha \mathop{\mathop{\sum' \nolimits}_{y = 1}^{p^\alpha}}_{g'(y) \equiv 0 (\bmod p^\alpha)} e\Bigl(\frac{g(y)}{p^{2\alpha}}\Bigr)
\end{equation}
where $g(y) = \frac{a y^5 + b}{y^3}$. Note $g'(y) = \frac{2 a y^7 - 3 b y^2}{y^6}$. Since $(y,p)=1$, $g'(y) \equiv 0 (\bmod p^\alpha)$ when
\begin{equation} \label{cond}
2 a y^5 - 3 b \equiv 0 (\bmod p^\alpha).
\end{equation}
If $(b,p) = p$, then $(a,p)=1$ and (\ref{cond}) has no solution with $(y,p) = 1$ unless $p=2$ in which case
\[
a y^5 - 3 (\frac{b}{2}) \equiv 0 (\bmod 2^{\alpha - 1})
\]
has at most one solution by Hensel's lemma. So we can suppose $(b,p)=1$. If $p=2$, (\ref{cond}) has no solution. If $p = 3$, then (\ref{cond}) has no solution unless $3 || a$ in which case
\[
2 (\frac{a}{3}) y^5 \equiv b (\bmod 3^{\alpha - 1})
\]
has at most five solutions by using primitive root $(\bmod 3^{\alpha - 1})$. If $p > 3$, then (\ref{cond}) has no solution unless $(a,p)=1$ in which case it has at most five solutions by using primitive root $(\bmod p^\alpha)$. Therefore we can conclude that
\begin{equation} \label{con-even}
|S(a,b;p^{2 \alpha})| \le 5 p^\alpha \text{ if } (a,b,p)=1.
\end{equation}

Now if $\beta = 2\alpha+1$ with $\alpha \ge 1$, by Lemma 12.3 in [\ref{IK}], we have
\begin{equation} \label{exp-oddpower}
S(a,b,p^{2\alpha+1}) = p^\alpha \mathop{\mathop{\sum' \nolimits}_{y = 1}^{p^\alpha}}_{g'(y) \equiv 0 (\bmod p^\alpha)} e\Bigl(\frac{g(y)}{p^{2\alpha+1}}\Bigr) G_p(y)
\end{equation}
where $g(y) = \frac{a y^5 + b}{y^3}$,
\[
G_p(y) = \sum_{z=1}^{p} e\Bigl(\frac{d(y) z^2 + g'(y) p^{-\alpha} z}{p}\Bigr)
\]
and $d(y) = \frac{g''(y)}{2}$. Note $g'(y) = \frac{2 a y^7 - 3 b y^2}{y^6}$ and $g''(y) = \frac{2 a y^8 + 12 b y^3}{y^8}$. So $d(y) = \frac{a y^8 + 6 b y^3}{y^8}$. If $(b,p)=p$, then by the analysis in the case $\beta = 2\alpha$, the sum in (\ref{exp-oddpower}) is empty unless $p=2$ in which case one has $|S(a,b,2^{2\alpha+1})| \le 2^{\alpha + 1}$. Now suppose $(b,p)=1$. If $p \neq 3, 5$, we consider those $y$ such that $(y,p)=1$ and satisfy (\ref{cond}). Then one can show that $p \not| 2d(y)$ (for otherwise $p | 2a y^5 + 12 b$ and $p | 2ay^5 - 3b$ imply $p | 15 b$, a contradiction.) Then the Gauss sum $G_p(y)$ satisfies $|G_p(y)| \le p^{1/2}$ (see equation (12.37) on page 322 of [\ref{IK}] for example) and hence $|S(a,b;p^{2\alpha + 1})| \le 5 p^{\alpha + 1/2}$ as there are at most five solutions to (\ref{cond}). If $p = 3$ or $5$, $|G_p(y)| \le 5$ and hence $|S(a,b;p^{2\alpha + 1})| \le 25 p^{\alpha + 1/2}$. In any case, we have
\begin{equation} \label{con-odd}
|S(a,b;p^{2\alpha + 1})| \le 25 p^{\alpha + 1/2} \text{ if } (a,b,p)=1.
\end{equation}

Hence combining (\ref{exp-prime}), (\ref{con-even}) and (\ref{con-odd}), we have
\[
|S(a,b,p^\beta)| \le 25 (a,b,p^\beta)^{1/2} p^{\beta/2}
\]
and, by (\ref{chain}), we have
\begin{lem} \label{lem3}
\[
|S(a,b,q)| \le 25^{\omega(q)} (a,b,q)^{1/2} q^{1/2} \ll_\epsilon (a,b,q)^{1/2} q^{1/2 + \epsilon}.
\]
\end{lem}

\section{Proof of Theorem \ref{firstthm}}

Recall $N_{A,B}(K,L) =: \# \{(a,b) : a^2 b^3 \equiv l (\bmod q), A < a \le A+K, B < b \le B+L \}$. Note that if $a^2 b^3 \equiv l (\bmod q)$, then $a^2 b^4 \equiv l b (\bmod q)$ and $b \equiv \overline{l} (a b^2)^2 (\bmod q)$. So
\[
N_{A,B}(K,L) = \# \{(a,b,w): b \equiv \overline{l} w^2 (\bmod q), a \equiv l \overline{w}^3 (\bmod q), A < a \le A+K, B < b \le B+L, 1 \le w \le q, (w,q)=1 \}.
\]

By orthogonality of additive characters,
\[
N_{A,B}(K,L) = \frac{1}{q^2} \sum_{m = 1}^{q} \sum_{n = 1}^{q} \sum_{A < a \le A+K} \sum_{B < b \le B+L} \mathop{\sum' \nolimits}_{w} e\Bigl(\frac{m (l^2 \overline{w}^3 - a)}{q}\Bigr) e\Bigl(\frac{n (\overline{l} w^2 - b)}{q}\Bigr).
\]
Separating the main contribution from $(m,n)=(q,q)$, we have
\begin{align*}
N_{A,B}(K,L) =& \frac{\phi(q) ([A+K]-[A]) ([B+L]-[B])}{q^2} \\
&+ \frac{K}{q^2} \sum_{n = 1}^{q-1} \sum_{B < b \le B+L} \mathop{\sum' \nolimits}_{w} e\Bigl(\frac{n (\overline{l} w^2 - b)}{q}\Bigr) + \frac{L}{q^2} \sum_{m = 1}^{q-1} \sum_{A < a \le A+K} \mathop{\sum' \nolimits}_{w} e\Bigl(\frac{m (l^2 \overline{w}^3 - a)}{q}\Bigr) \\
&+ \frac{1}{q^2} \sum_{m = 1}^{q-1} \sum_{n = 1}^{q-1} \sum_{A < a \le A+K} \sum_{B < b \le B+L} \mathop{\sum' \nolimits}_{w} e\Bigl(\frac{m (l \overline{w}^3 - a)}{q}\Bigr) e\Bigl(\frac{n (\overline{l} w^2 - b)}{q}\Bigr) := S_1 + S_2 + S_3 + S_4.
\end{align*}
Firstly,
\[
S_2 = \frac{K}{q^2} \sum_{n = 1}^{q-1} \sum_{B < b \le B+L} e\Bigl(-\frac{n b}{q}\Bigr) S(\overline{l} n, 0; q) \ll_\epsilon \frac{K}{q^{1/2 - \epsilon}}
\]
by Lemma \ref{lem1} and Lemma \ref{lem3}. Similarly,
\[
S_3 \ll_\epsilon \frac{L}{q^{1/2 - \epsilon}}.
\]
Finally,
\begin{align*}
S_4 =& \frac{1}{q^2} \sum_{m = 1}^{q-1} \sum_{A < a \le A+K} e\Bigl(-\frac{m a}{q}\Bigr) \sum_{n = 1}^{q-1} \sum_{B < b \le B+L} e\Bigl(-\frac{n b}{q}\Bigr) S(m l^2, n \overline{l}; q) \\
\ll_\epsilon& \frac{1}{q^2} \sum_{m = 1}^{q-1} \Big| \sum_{A < a \le A+K} e\Bigl(-\frac{m a}{q}\Bigr) \Big| \sum_{n = 1}^{q-1} \Big| \sum_{B < b \le B+L} e\Bigl(-\frac{n b}{q}\Bigr) \Big| (m, n, q)^{1/2} q^{1/2 + \epsilon} \\
\ll_\epsilon& \frac{1}{q^{3/2-\epsilon}} \sum_{m = 1}^{q-1} \Big| \sum_{A < a \le A+K} e\Bigl(-\frac{m a}{q}\Bigr) \Big| \sum_{n = 1}^{q-1} \Big| \sum_{B < b \le B+L} e\Bigl(-\frac{n b}{q}\Bigr) \Big| (n, q)^{1/2} \ll_\epsilon q^{1/2+\epsilon}
\end{align*}
by Lemmas \ref{lem1}, \ref{lem2} and \ref{lem3}. Therefore
\[
N_{A,B}(K,L) = \frac{\phi(q) K L}{q^2} + O\Bigl(\frac{K+L}{q^{1/2-\epsilon}} + q^{1/2+\epsilon}\Bigr)
\]
which gives Theorem \ref{firstthm}.
\section{Proof of Theorem \ref{interthm}}

Suppose $\lambda, \mu < 1/5$.
\begin{align*}
S_q(X; l) =& \mathop{\sum_{a^2 b^3 \le X}}_{a^2 b^3 \equiv l \; (\bmod q)} 1 \\
=& \sum_{b \le X^{\mu}} \mathop{\sum_{a \le X^{1/2} /
b^{3/2}}}_{a^2 b^3 \equiv l (\bmod q)} 1 + \sum_{a \le X^{\lambda}}
\mathop{\sum_{b \le X^{1/3} / a^{2/3}}}_{a^2 b^3 \equiv l (\bmod q)}
1 - \sum_{b \le X^{\mu}} \mathop{\sum_{a \le X^{\lambda}}}_{a^2 b^3
\equiv l (\bmod q)} 1 + \mathop{\sum_{X^\lambda < a} \sum_{X^\mu < b}}_{a^2 b^3 \le X, a^2 b^3 \equiv l (\bmod q)} 1 \\
=:& T_1 + T_2 - T_3 + T_4.
\end{align*}

\begin{align*}
T_3 =& \mathop{\sum_{b \le X^{\mu}}}_{(b,q) = 1} \mathop{\sum_{a \le
X^{\lambda}}}_{a^2 \equiv l \overline{b}^3 (\bmod q)} 1 = \sum_{b \le
X^{\mu}} N_2(l b; q) \Bigl[ \frac{X^{\lambda}}{q} + O(1) \Bigr] \\
=& \frac{X^{\lambda}}{q} \sum_{b \le X^{\mu}} N_2(l b; q) + O_\epsilon \Bigl(\sum_{b \le
X^{\mu}} N_2(l b; q)\Bigr) = \frac{\phi(q)}{q^2} X^{\lambda + \mu} + O_\epsilon\Bigl(\frac{X^\lambda}{q^{1/2 - \epsilon}} \Bigr) + O_\epsilon( X^{\mu} q^\epsilon)
\end{align*}
by definition of $N_2(n;q)$, Lemma \ref{lem2} and $N_2(n; q) \ll_\epsilon q^{\epsilon}$ (See Nagell [\ref{N}]).
Secondly,
\begin{align*}
T_1 =& \mathop{\sum_{b \le X^{\mu}}}_{(b,q) = 1} \mathop{\sum_{a \le
X^{1/2} / b^{3/2}}}_{a^2 \equiv l \overline{b}^3 (\bmod q)} 1 =
\sum_{b \le X^{\mu}} N_2(l b; q) \Bigl[
\frac{X^{1/2}}{b^{3/2} q} + O(1) \Bigr]\\
=& \frac{X^{1/2}}{q} \sum_{b \le X^{\mu}} \frac{N_2(l b; q)}{b^{3/2}}
+ O_\epsilon (X^{\mu} q^\epsilon) \\
=& \frac{X^{1/2}}{q} \sum_{b = 1}^{\infty} \frac{N_2(l b; q)}{b^{3/2}} - \frac{X^{1/2}}{q} \int_{X^{\mu}}^{\infty} \frac{1}{u^{3/2}} d \sum_{b \le u} N_2(l b;q) + O_\epsilon (X^{\mu} q^\epsilon) \\
=& C_q(l) \frac{X^{1/2}}{q} - 2 \frac{\phi(q) X^{1/2 - \mu/2}}{q^2} - \frac{X^{1/2}}{q} \int_{X^{\mu}}^{\infty} \frac{1}{u^{3/2}} d \Bigl(\sum_{b \le u} N_2(l b;q) - \frac{\phi(q)}{q} u\Bigr) + O_\epsilon (X^{\mu}
q^\epsilon) \\
=& C_q(l) \frac{X^{1/2}}{q} - 2 \frac{\phi(q) X^{1/2 - \mu/2}}{q^2} - \frac{3}{2} \frac{X^{1/2}}{q} \int_{X^{\mu}}^{\infty} \frac{\sum_{b \le u} N_2(l b;q) - \frac{\phi(q)}{q} u}{u^{5/2}} du + O_\epsilon (X^{\mu}
q^\epsilon) \\
=& C_q(l) \frac{X^{1/2}}{q} - 2 \frac{\phi(q)}{q^2} X^{1/2 - \mu/2} + O\Bigl(\frac{X^{1/2 - 3\mu/2}}{q^{1/2 - \epsilon}} \Bigr) + O_\epsilon (X^{\mu} q^\epsilon)
\end{align*}
by definition of $N_2(n;q)$, Lemma \ref{lem2} and $N_2(n; q) \ll_\epsilon q^{\epsilon}$.

Thirdly, let
\[
F(u) := \sum_{l_1 \le u} N_3(l l_1; q) - \frac{\phi(q)}{q} u.
\]
We can have
\begin{align*}
T_2 =& \mathop{\sum_{a \le X^{\lambda}}}_{(a,q) = 1} \mathop{\sum_{b \le
X^{1/3} / a^{2/3}}}_{b^3 \equiv l \overline{a}^2 (\bmod q)} 1 =
\sum_{a \le X^{\lambda}} N_3(l a; q) \Bigl[ \frac{X^{1/3}}{a^{2/3} q} + O(1) \Bigr] \\
=& \frac{X^{1/3}}{q} \sum_{a \le X^{\lambda}} \frac{N_3(l a;
q)}{a^{2/3}} + O\Bigl( \sum_{a \le X^{\lambda}} N_3(l a; q) \Bigr) \\
=& \frac{X^{1/3}}{q} \int_{1}^{X^{\lambda}} \frac{1}{u^{2/3}} d F(u) + \frac{X^{1/3}}{q} \frac{\phi(q)}{q} \int_{1}^{X^{\lambda}} \frac{1}{u^{2/3}} du + O_\epsilon (X^{\lambda} q^\epsilon) \\
=& \frac{X^{1/3}}{q} \Bigl[ O\Bigl(\frac{q^{1/2+\epsilon}}{(X^{\lambda})^{2/3}} \Bigr) + \frac{\phi(q)}{q} + \frac{2}{3}
\int_{1}^{X^{\lambda}} \frac{F(u)}{u^{5/3}} du \Bigr] + \frac{\phi(q) X^{1/3}}{q^2} [3 X^{\lambda/3} - 3] + O_\epsilon
(X^{\lambda} q^\epsilon) \\
=& 3 \frac{\phi(q)}{q^2} X^{1/3 + \lambda/3} + \Bigl[ -2 \frac{\phi(q)}{q} + \frac{2}{3}
\int_{1}^{X^{\lambda}} \frac{F(u)}{u^{5/3}} du \Bigr] \frac{X^{1/3}}{q} + O\Bigl(\frac{X^{1/3 - 2\lambda/3}}{q^{1/2 - \epsilon}}\Bigr) + O_\epsilon(X^{\lambda} q^\epsilon) \\
=& 3 \frac{\phi(q)}{q^2} X^{1/3 + \lambda/3} + \Bigl[ -2 \frac{\phi(q)}{q} + \frac{2}{3}
\int_{1}^{\infty} \frac{F(u)}{u^{5/3}} du \Bigr] \frac{X^{1/3}}{q} + O\Bigl(\frac{X^{1/3 - 2\lambda/3}}{q^{1/2 - \epsilon}}\Bigr) + O_\epsilon(X^{\lambda} q^\epsilon)
\end{align*}
by definition of $N_3(n;q)$, Lemma \ref{lem2} and $N_3(n;q) \ll_\epsilon q^\epsilon$.

\bigskip

Finally, we deal with $T_4$. We are going to partition the region $\{(a,b): X^\lambda < a, X^\mu < b, a^2 b^3 \le X \}$ into rectangles with some left-over regions. We follow closely the construction in [\ref{MV}]. Here we require that the choices for $\lambda$ and $\mu$ are such that $\log_2 X^{1/2 - 3\mu/2 - \lambda}$ is an integer.

Firstly, we begin with the rectangles
\[
R_{i} =: (2^{i-1} X^\lambda, 2^i X^\lambda] \times \Bigl(X^\mu, (\frac{X}{(2^i X^\lambda)^2})^{1/3}\Bigr] \; \; \text{ for } 1 \le i < I = \log_2 X^{1/2 - 3\mu/2 - \lambda}.
\]
In the remaining regions
\[
S_{i} =: \Bigl\{(a,b): 2^{i-1} X^\lambda < a \le 2^i X^\lambda, b > (\frac{X}{(2^i X^\lambda)^2})^{1/3}, a^2 b^3 \le X \Bigr\} \text{ for } i = 1,2, ... , I;
\]
we place additional rectangles $R_{ijk}$. In $S_{i}$ we put
\[
R_{i11} = \Bigl(2^{i-1} X^\lambda, 2^{i-1} (1 + \frac{1}{2}) X^\lambda\Bigr] \times \Bigl((\frac{X}{(2^i X^\lambda)^2})^{1/3}, (\frac{X}{(2^{i-1} (1 + \frac{1}{2}) X^\lambda)^2})^{1/3}\Bigr].
\]
Left over in $S_i$ are two regions of the same kind into each of which we place a further rectangle and so on. Thus on the $j-$th occasion we place $2^{j-1}$ rectangles $R_{ijk}$ ($1 \le k \le 2^{j-1}$) where
\[
R_{ijk} = \Bigl(2^{i-1}(1 + \frac{2k-2}{2^j}) X^\lambda, 2^{i-1}(1 + \frac{2k-1}{2^j}) X^\lambda\Bigr] \times \Bigl((\frac{X}{(2^{i-1} (1 + \frac{2k}{2^j}) X^\lambda)^2})^{1/3}, (\frac{X}{(2^{i-1} (1 + \frac{2k-1}{2^j}) X^\lambda)^2})^{1/3}\Bigr].
\]
The left over regions are
\[
S_{ijk} = \Bigl\{(a,b) : 2^{i-1}(1 + \frac{2k-2}{2^j}) X^\lambda < a \le 2^{i-1}(1 + \frac{2k-1}{2^j}) X^\lambda, b > \Bigl(\frac{X}{(2^{i-1}(1 + \frac{2k-1}{2^j}) X^\lambda)^2}\Bigr)^{1/3}, a^2 b^3 \le X \Bigr\}
\]
and
\[
S_{ijk}' = \Bigl\{(a,b) : 2^{i-1}(1 + \frac{2k-1}{2^j}) X^\lambda < a \le 2^{i-1}(1 + \frac{2k}{2^j}) X^\lambda, b > \Bigl(\frac{X}{(2^{i-1}(1 + \frac{2k}{2^j}) X^\lambda)^2}\Bigr)^{1/3}, a^2 b^3 \le X \Bigr\}.
\]
We keep doing this until $j = J$. Hence, by Theorem \ref{firstthm},
\begin{align*}
T_4 =& \sum_{i = 1}^{I-1} \mathop{\sum_{(a,b) \in R_i}}_{a^2 b^3 \equiv l (\bmod q)} 1 + \sum_{j = 1}^{J} \sum_{i=1}^{I} \sum_{k = 1}^{2^{j-1}} \mathop{\sum_{(a,b) \in R_{ijk}}}_{a^2 b^3 \equiv l (\bmod q)} 1 + \sum_{i=1}^{I} \sum_{k = 1}^{2^{J-1}} \mathop{\sum_{(a,b) \in S_{iJk}}}_{a^2 b^3 \equiv l (\bmod q)} 1 + \sum_{i=1}^{I} \sum_{k = 1}^{2^{J-1}} \mathop{\sum_{(a,b) \in S_{iJk}'}}_{a^2 b^3 \equiv l (\bmod q)} 1 \\
=& \frac{\phi(q)}{q^2} \text{Area of } \{(a,b) : X^\lambda < a, X^\mu < b, a^2 b^3 \le X\} + O_\epsilon(I 2^J q^{1/2+\epsilon}) + O\Bigl(\sum_{i=1}^{I} \sum_{k = 1}^{2^{J-1}} \frac{\text{Area of } S_{iJk}}{q}\Bigr) \\
&+ O\Bigl(\sum_{i=1}^{I} \sum_{k = 1}^{2^{J-1}} \frac{\text{Area of }S_{iJk}'}{q}\Bigr) + O\Bigl(\sum_{i=1}^{I} \sum_{k = 1}^{2^{J-1}} \mathop{\sum_{(a,b) \in S_{iJk}}}_{a^2 b^3 \equiv l (\bmod q)} 1\Bigr) + O\Bigl(\sum_{i=1}^{I} \sum_{k = 1}^{2^{J-1}} \mathop{\sum_{(a,b) \in S_{iJk}'}}_{a^2 b^3 \equiv l (\bmod q)} 1\Bigr)
\end{align*}
Note that
\[
S_{iJk} \subset \Bigl(2^{i-1}(1 + \frac{2k-2}{2^J}) X^\lambda, 2^{i-1}(1 + \frac{2k-1}{2^J}) X^\lambda\Bigr] \times \Bigl((\frac{X}{(2^{i-1}(1 + \frac{2k-1}{2^J}) X^\lambda)^2})^{1/3}, (\frac{X}{(2^{i-1}(1 + \frac{2k-2}{2^J}) X^\lambda)^2})^{1/3}\Bigr].
\]
So the area of $S_{iJk}$ is at most
\[
\frac{2^{i-1} X^\lambda}{2^J} \Bigl[\Bigl(\frac{X}{(2^{i-1}(1 + \frac{2k-2}{2^J}) X^\lambda)^2}\Bigr)^{1/3} - \Bigl(\frac{X}{(2^{i-1}(1 + \frac{2k-1}{2^J}) X^\lambda)^2}\Bigr)^{1/3}\Bigr] \ll \frac{2^{i-1} X^\lambda}{2^J} \cdot \frac{X^{1/3}}{(2^{i-1} X^\lambda)^{2/3}} \cdot \frac{1}{2^J}.
\]
So the second and fourth error terms are
\[
\ll_\epsilon \sum_{i=1}^{I} \sum_{k = 1}^{2^{J-1}} \frac{(2^{i-1} X^\lambda)^{1/3} X^{1/3}}{q} \frac{1}{2^{2J}} + q^{1/2 + \epsilon} \ll \frac{X^{1/2 - \mu/2}}{2^J q} + 2^J q^{1/2 + \epsilon} \log X.
\]
Similarly, we have the same error bound for the third and the fifth error terms. Therefore
\begin{align*}
T_4 =& \frac{\phi(q)}{q^2} \int_{X^\lambda}^{X^{1/2-3\mu/2}} \Bigl(\frac{X}{a^2}\Bigr)^{1/3} - X^\mu da + O_\epsilon \Bigl(\frac{X^{1/2 - \mu/2}}{2^J q} + 2^J q^{1/2 + \epsilon} \log X \Bigr) \\
=& \frac{\phi(q)}{q^2} \Bigl[2X^{1/2 - \mu/2} - 3X^{1/3 + \lambda/3} - X^{\lambda+\mu}\Bigr] + O_\epsilon \Bigl(\frac{X^{1/2 - \mu/2}}{2^J q} + 2^J q^{1/2 + \epsilon} \log X \Bigr).
\end{align*}
Hence
\[
S_q(X;l) = C_q(l) \frac{X^{1/2}}{q} + D_q(l) \frac{X^{1/3}}{q} + O_\epsilon \Bigl(\frac{X^{1/2 - 3\mu/2}}{q^{1/2 - \epsilon}} + \frac{X^{1/3 - 2\lambda/3}}{q^{1/2 - \epsilon}} + X^\lambda q^\epsilon + X^\mu q^\epsilon + \frac{X^{1/2 - \mu/2}}{2^J q} + 2^J q^{1/2 + \epsilon} \log X \Bigr).
\]
We are going to pick $X^\lambda$ to have the same size as $X^\mu$ with $\mu \le 1/5$. Then the first error term dominates the second one. We will also pick $2^J$ such that $\frac{X^\mu}{2^J q^{1/2}} \ll 1$. Then the first error term dominates the fifth one. So
\[
S_q(X;l) = C_q(l) \frac{X^{1/2}}{q} + D_q(l) \frac{X^{1/3}}{q} + O_\epsilon \Bigl(\frac{X^{1/2 - 3\mu/2}}{q^{1/2 - \epsilon}} + X^\mu q^\epsilon + 2^J q^{1/2 + \epsilon} \Bigr).
\]
If $q \le X^{2\mu}$, then we can pick $2^J$ of size $\frac{X^\mu}{q^{1/2}}$ and bigger than $1$ and obtain
\[
S_q(X;l) = C_q(l) \frac{X^{1/2}}{q} + D_q(l) \frac{X^{1/3}}{q} + O_\epsilon \Bigl(\frac{X^{1/2 - 3\mu/2}}{q^{1/2 - \epsilon}} + X^\mu q^\epsilon \Bigr).
\]
If $q > X^{2\mu}$, we simply pick $2^J = 2$ and get
\[
S_q(X;l) = C_q(l) \frac{X^{1/2}}{q} + D_q(l) \frac{X^{1/3}}{q} + O_\epsilon \Bigl(\frac{X^{1/2 - 3\mu/2}}{q^{1/2 - \epsilon}} + q^{1/2 + \epsilon} \Bigr).
\]
The above together gives Theorem \ref{interthm}.
\section{Proof of Theorem \ref{finalthm}}

First recall that any squarefull number
$n$ can be written uniquely as $n = a^2 b^3$ with $b$ is squarefree.
Hence, by Theorem \ref{oldthm1}, $C_q(l a^2) = C_q(l)$ for any $(a,q) =
1$, and $C_q(l) \ll_\epsilon q^\epsilon$ as $N_2(n; q) \ll_\epsilon
q^{\epsilon}$,
\begin{align*}
T_q(x;l) :=& \mathop{\mathop{\sum_{a^2 b^3 \le x}}_{b \text{
squarefree}}}_{a^2 b^3 \equiv l (\bmod q)} 1 = \mathop{\sum_{a^2 b^3
\le x}}_{a^2 b^3 \equiv l (\bmod q)} \sum_{d^2 | b} \mu(d) =
\mathop{\sum_{a^2 d^6 b'^3 \le x}}_{a^2 d^6 b'^3
\equiv l (\bmod q)} \mu(d) \\
=& \mathop{\sum_{d \le x^{1/6}/q^{1/(12 \mu)}}}_{(d,q) = 1} \mu(d)
\mathop{\sum_{a^2 b'^3 \le x/d^6}}_{a^2 b'^3 \equiv l \overline{d}^6 (\bmod q)} 1 +
\mathop{\sum_{x^{1/6}/q^{1/(12 \mu)} < d \le x^{1/6}/q^{\nu}}}_{(d,q) = 1} \mu(d)
\mathop{\sum_{a^2 b'^3 \le x/d^6}}_{a^2 b'^3 \equiv l \overline{d}^6 (\bmod q)} 1 \\
&+ \mathop{\sum_{x^{1/6}/q^{\nu} < d \le x^{1/6}}}_{(d,q) = 1} \mu(d)
\mathop{\sum_{a^2 b'^3 \le x/d^6}}_{a^2 b'^3 \equiv l \overline{d}^6 (\bmod q)} 1 \\
=:& U_1 + U_2 + U_3
\end{align*}
For $U_1$, one can check that the condition on $d$ implies $q \le (\frac{x}{d^6})^{2\mu}$. So we can apply Theorem \ref{interthm} and get
\begin{align*}
U_1 =& \mathop{\sum_{d \le x^{1/6}/q^{1/(12\mu)}}}_{(d,q) = 1} \mu(d)
\Bigl[ C_q(l) \frac{x^{1/2}}{d^3 q} + D_q(l) \frac{x^{1/3}}{d^2 q} + O_\epsilon \Bigl(\frac{x^{1/2 - 3\mu/2}}{d^{3 - 9\mu} q^{1/2 - \epsilon}} + \frac{x^\mu q^\epsilon}{d^{6 \mu}} \Bigr) \Bigr] \\
=& C_q(l) \mathop{\sum_{d = x^{1/6}/q^{1/(12\mu)}}}_{(d,q) = 1} \frac{\mu(d)}{d^3} \frac{x^{1/2}}{q} + D_q(l) \mathop{\sum_{d \le x^{1/6}/q^{1/(12\mu)}}}_{(d,q) = 1} \frac{\mu(d)}{d^2} \frac{x^{1/3}}{q} + O_\epsilon \Bigl(\frac{x^{1/2 - 3\mu/2}}{q^{1/2 - \epsilon}} + x^\mu q^\epsilon \Bigr)
\end{align*}
as $C_q(l) \ll q^\epsilon$ and $D_q(l) \ll q^{\epsilon}$. By Theorem \ref{interthm},
\begin{align*}
U_2 =& \mathop{\sum_{x^{1/6}/q^{1/(12 \mu)} < d \le x^{1/6}/q^{\nu}}}_{(d,q) = 1} \mu(d)
\mathop{\sum_{a^2 b'^3 \le x/d^6}}_{a^2 b'^3 \equiv l \overline{d}^6 (\bmod q)} 1 \\
=& \mathop{\sum_{x^{1/6}/q^{1/(12 \mu)} < d \le x^{1/6}/q^{\nu}}}_{(d,q) = 1} \mu(d) \Bigl[ C_q(l) \frac{x^{1/2}}{d^3 q}  + D_q(l) \frac{x^{1/3}}{d^2 q} + O_\epsilon \Bigl( \Bigl(\frac{x^{1/2 - 3\mu/2}}{d^{3 - 9\mu} q^{1/2}} + q^{1/2} \Bigr) q^\epsilon \Bigr) \Bigr]\\
=& C_q(l) \mathop{\sum_{x^{1/6}/q^{1/(12 \mu)} < d \le x^{1/6}/q^{\nu}}}_{(d,q) = 1} \frac{\mu(d)}{d^3} \frac{x^{1/2}}{q} + D_q(l) \mathop{\sum_{x^{1/6}/q^{1/(12 \mu)} < d \le x^{1/6}/q^{\nu}}}_{(d,q) = 1} \frac{\mu(d)}{d^2} \frac{x^{1/3}}{q}\\
&+ O_\epsilon \Bigl( \Bigl(\frac{x^{1/6}}{q^{1/2 - (2 - 9\mu)/(12\mu)}} + x^{1/6} q^{1/2 - \nu} \Bigr) q^\epsilon \Bigr).
\end{align*}
Finally, by Theorem \ref{oldthm1},
\begin{align*}
U_3 \le& \mathop{\sum_{x^{1/6}/q^{\nu} < d \le x^{1/6}}}_{(d,q) = 1}
\mathop{\sum_{a^2 b'^3 \le x/d^6}}_{a^2 b'^3 \equiv l \overline{d}^6 (\bmod q)} 1 \\
\ll_\epsilon& \mathop{\sum_{x^{1/6}/q^{\nu} < d \le x^{1/6}}}_{(d,q) = 1} \Bigl[ \frac{(x/d^6)^{1/2}}{q} + \Bigl(\frac{x}{d^6}\Bigr)^{1/5} \Bigr] q^\epsilon \\
\ll& (x^{1/6} q^{2\nu - 1} + x^{1/6} q^{\nu/5}) q^\epsilon.
\end{align*}
To minimize the errors, we set $1/2 - \nu = \nu/5$ which means $\nu = 5/12$; and set $x^{1/2 - 3\mu/2} / q^{1/2} = x^\mu$ which means $\mu = 1/5 - \theta/5$ where $q = x^\theta$. Consequently, we have
\[
T_q(x;l) = C_q(l) \mathop{\sum_{d = 1}^{\infty}}_{(d,q) = 1} \frac{\mu(d)}{d^3} \frac{x^{1/2}}{q} + D_q(l) \mathop{\sum_{d = 1}^{\infty}}_{(d,q) = 1} \frac{\mu(d)}{d^2} \frac{x^{1/3}}{q} + O_\epsilon \Bigl( \Bigl( x^{1/6} q^{1/12} + \frac{x^{1/5}}{q^{1/5}} \Bigr) q^\epsilon \Bigr)
\]
which is Theorem \ref{finalthm}.

Tsz Ho Chan \\
Department of Arts and Sciences \\
Victory University \\
255 N. Highland Street \\
Memphis, TN 38111 \\
U.S.A. \\
White Station High School \\
514 S. Perkins Road \\
Memphis, TN 38117 \\
U.S.A. \\
thchan6174@gmail.com \\

\end{document}